\def\yes{\if00}
\def\no{\if01}
\def\iftwelvept{\yes}
\def\ifusepdf{\no}
\def\ifpsfont{\yes}

\iftwelvept
\documentclass[leqno,12pt]{amsart}
\else
\documentclass[leqno]{amsart}
\fi
\usepackage{amssymb}
\usepackage{amscd}
\usepackage{latexsym}
\usepackage{verbatim}

\ifusepdf
\usepackage{hyperref}
\else\fi
\ifpsfont
\usepackage[T1]{fontenc}
\usepackage{times}
\else\fi

\iftwelvept
\else\fi

\theoremstyle{plain}
\newtheorem{Theorem}{Theorem}[section]
\newtheorem{Proposition}[Theorem]{Proposition}
\newtheorem{Lemma}[Theorem]{Lemma}
\newtheorem{Corollary}[Theorem]{Corollary}
\newtheorem{Claim}{Claim}[Theorem]

\newtheorem{Sublemma}[Claim]{Sublemma}

\theoremstyle{definition}

\renewcommand{\theTheorem}{\arabic{section}.\arabic{Theorem}}
\renewcommand{\theClaim}{\arabic{section}.\arabic{Theorem}.\arabic{Claim}}
\renewcommand{\theequation}{\arabic{section}.\arabic{Theorem}.\arabic{Claim}}

\def\rom{\textup}
\newcommand{\ZZ}{{\mathbb{Z}}}
\newcommand{\QQ}{{\mathbb{Q}}}
\newcommand{\RR}{{\mathbb{R}}}
\newcommand{\CC}{{\mathbb{C}}}
\newcommand{\PP}{{\mathbb{P}}}

\newcommand{\OO}{{\mathcal{O}}}
\newcommand{\MM}{\operatorname{\mathit{M}}}
\newcommand{\bMM}{\operatorname{\bar{\mathit{M}}}}

\newcommand{\bcMM}{\operatorname{\bar{\mathcal{M}}}}
\newcommand{\Ker}{\operatorname{Ker}}

\newcommand{\cl}{\operatorname{cl}}
\newcommand{\Image}{\operatorname{Im}}
\newcommand{\Nm}{\operatorname{Nm}}
\newcommand{\codim}{\operatorname{codim}}

\newcommand{\End}{\operatorname{End}}
\newcommand{\Pic}{\operatorname{Pic}}
\newcommand{\NS}{\operatorname{NS}}
\newcommand{\NSnum}{\operatorname{NS}^{\nu}}
\newcommand{\Cycle}{Z}

\newcommand{\Spec}{\operatorname{Spec}}
\newcommand{\Chow}{\operatorname{CH}}

\newcommand{\Proof}{{\sl Proof.}\quad}
\newcommand{\QED}{{\unskip\nobreak\hfil\penalty50\quad\null\nobreak\hfil
{$\Box$}\parfillskip0pt\finalhyphendemerits0\par\medskip}}
\newcommand{\rest}[2]{\left.{#1}\right\vert_{{#2}}}

\begin{document}

\title[The $\QQ$-Picard group of the moduli space of curves]
{The $\QQ$-Picard group of the moduli space of curves \\
in positive characteristic}
\author{Atsushi Moriwaki}
\address{Department of Mathematics, Faculty of Science,
Kyoto University, Kyoto, 606-8502, Japan}
\email{moriwaki@kusm.kyoto-u.ac.jp}
\date{15/Dec/2000, 11:45PM(JP), (Version 2.0)}
\begin{abstract}
In this note, we prove that the $\QQ$-Picard group of the moduli
space of $n$-pointed stable curves of genus $g$ over an algebraically
closed field is generated by the tautological classes.
We also prove that the cycle map to the 2nd \'etale
cohomology group is bijective.
\end{abstract}


\maketitle


\section*{Introduction}
Let $k$ be an algebraically closed field,
$g$ and $n$ non-negative integers with $2g-2+n > 0$, and
$\bMM_{g,n}$ (resp. $\MM_{g,n}$) the moduli space of $n$-pointed
stable (resp. smooth) curves of genus $g$ over $k$.
We denote by $\Pic(\bMM_{g,n})_{\QQ}$ the $\QQ$-Picard group of $\bMM_{g,n}$;
that is, $\Pic(\bMM_{g,n})_{\QQ} = \Pic(\bMM_{g,n}) \otimes \QQ$.
Let $\lambda$ be the Hodge class,
$\psi_1, \ldots, \psi_n$ the classes of
$\QQ$-line bundles
given by the pull-back of 
the relative dualizing sheaf of the universal
curve over $\bMM_{g,n}$ in terms of $n$ sections, and
$\{ \delta_{t} \}_{t \in T}$ the boundary classes in
$\Pic(\bMM_{g,n})_\QQ$
(for details, see \S\S\ref{subsec:tautological:classes}).
These classes
\[
\text{$\lambda, \psi_1, \ldots, \psi_n$ and
$\delta_{t}$'s ($t \in T$)}
\]
are called {\em the tautological classes of $\Pic(\bMM_{g,n})_\QQ$}.
It is well known (due to Harer) that
$\Pic(\bMM_{g,n})_\QQ$ is generated by the tautological classes
if the characteristic of $k$ is zero
(see Arbarello-Cornalba~\cite{AC}
for its proof by means of algebraic geometry).
In this note, we would like to show that
this still holds even if
the characteristic of $k$ is positive.
Namely, we have the following.

\renewcommand{\theTheorem}{}
\begin{Theorem}[cf. Theorem~\ref{thm:generators:NS:Mgn}]
\label{thm:generators:NS:Mgn:intro}
$\Pic(\bMM_{g,n})_\QQ$ is generated by the tautological classes
\[
\text{$\lambda, \psi_1, \ldots, \psi_n$ and
$\delta_{t}$'s \rom{(}$t \in T$\rom{)}}
\]
for any algebraically closed field $k$.
Moreover, the cycle map
\[
\Pic(\bMM_{g,n}) \otimes \QQ_{\ell} \to H_{et}^2(\bMM_{g,n}, \QQ_{\ell})
\]
is bijective for every prime $\ell$ invertible in $k$.
\end{Theorem}
\renewcommand{\theTheorem}{\arabic{section}.\arabic{Theorem}}

We prove the above theorem by using modulo $p$ reduction.
The outline of the proof is as follows:
Let $\bcMM_{g,n}$ be the algebraic stack classifying
$n$-pointed stable curves of genus $g$.
We compare the \'etale cohomology group
$H_{et}^2(\bcMM_{g,n}(k), \QQ_{\ell})$
over $k$ with the singular cohomology group $H^2(\bcMM_{g,n}(\CC), \QQ)$
of the analytic space $\bcMM_{g,n}(\CC)$
via a smooth Galois covering of
$\bcMM_{g,n}$ in terms of
a Teichm\"{u}ller level structure
due to Looijenga-Pikaart-de Jong-Boggi; namely,
\[
\dim_{\QQ_{\ell}} H_{et}^2(\bcMM_{g,n}(k), \QQ_{\ell})
\leq \dim_{\QQ} H^2(\bcMM_{g,n}(\CC), \QQ).
\]
Moreover, using the simple connectedness of
the moduli space of curves with a level $m \geq 3$
(due to Boggi-Pikaart),
we see that the cycle map
\[
\Pic(\bcMM_{g,n}(k)) \otimes \QQ_{\ell}
\to H_{et}^2(\bcMM_{g,n}(k), \QQ_{\ell})
\]
is injective. In this way, we obtain our theorem together with
the linear independence of the tautological classes.

After writing this note, Prof. Keel informed me
that Prof. de Jong had the same idea for the proof
of the above theorem.

\section{Notations and conventions}
\subsection{}
For a finite set $S$, we denote the number of it by $\vert S \vert$.

\subsection{}
If a group $G$ acts on a certain kind
of a mathematical object, the induced automorphism
by $g \in G$ is denoted by $[g]$.

\subsection{}
\label{subsec:num:NS}
Let $X$ be a proper algebraic space over an algebraically
closed field $k$.
Let $L_1$ and $L_2$ be line bundles on $X$.
We say $L_1$ is {\em numerically equivalent} to $L_2$,
denoted by $L_1 \equiv L_2$, if
$(L_1 \cdot C) = (L_2 \cdot C)$ for all
curves $C$ on $X$.
The group $\Pic(X)$ modulo
the numerical equivalence, denoted by $\NSnum(X)$,
is called {\em the numerical N\'eron-Severi group of $X$}.
Moreover, we denote $\NSnum(X) \otimes \QQ$ by
$\NSnum(X)_{\QQ}$.

\subsection{}
Let $X$ be an algebraic scheme over an algebraically
closed field $k$.
We define the {\em $\QQ$-Picard group $\Pic(X)_{\QQ}$ of $X$} to be
$\Pic(X)_{\QQ} = \Pic(X) \otimes \QQ$.
For $L_1, L_2 \in \Pic(X)$,
we say $L_1$ is {\em algebraically equivalent} to $L_2$,
denoted by $L_1 \sim_{alg} L_2$, if
there are a connected and smooth algebraic scheme $T$ over $k$, a line
bundle $\mathcal{L}$ on $T \times X$ and
$t_1, t_2 \in T$ such that
$\rest{\mathcal{L}}{\{t_1\} \times X} \simeq L_1$ and
$\rest{\mathcal{L}}{\{t_2\} \times X} \simeq L_2$.
The N\'eron-Severi group $\NS(X)$ of $X$ is defined by
$\Pic(X)$ modulo the algebraic equivalence.
In other words, $\NS(X) = \Pic_X(k)/(\Pic_X^0)_{red}(k)$,
where $\Pic_X$ is the Picard
scheme of $X$ and
$\Pic_X^0$ is the connected component containing $0$.
The group $\NS(X) \otimes \QQ$, denoted by $\NS(X)_\QQ$,
is called {\em the $\QQ$-N\'eron-Severi group of $X$}.
We assume that $X$ is projective over $k$.
It is not difficult to see that $L_1 \sim_{alg} L_2$ implies
$L_1 \equiv L_2$, so that
we have the natural surjective homomorphism
$\NS(X) \to \NSnum(X)$. It is well known (due to Matsusaka)
that the kernel of $\NS(X) \to \NSnum(X)$ is a finite
group. Thus, we can identify $\NS(X)_{\QQ}$ with
$\NSnum(X)_{\QQ}$.

\subsection{}
\label{subsec:alg:stack}
In this note, an algebraic stack always means a separated
algebraic stack over a locally noetherian scheme in the sense
of Deligne-Mumford.
Let $X$ be a algebraic stack over a locally noetherian scheme $S$.
For an algebraically closed field $L$ and
a morphism $\Spec(L) \to S$, the coarse moduli space
of $X \times_{S} \Spec(L)$ is denoted by
$X_L$ (cf. \cite[Chapter~I, Theorem~4.10]{FC} and \cite[Corollary~1.3]{KM}).

\subsection{}
\label{subsec:tautological:classes}
Let $g$ and $n$ be non-negative integers with $2g-2+n > 0$, and
$\bMM_{g,n}$ (resp. $\MM_{g,n}$) the moduli space of $n$-pointed
stable (resp. smooth) curves of genus $g$ over
an algebraically closed field.  Roughly speaking,
the $\QQ$-line bundles $\lambda$ and $\psi_1, \ldots, \psi_n$ on
$\bMM_{g,n}$ are defined as follows:
Let $\pi : \bMM_{g,n+1} \to \bMM_{g,n}$ be the universal curve
of $\bMM_{g,n}$, and $s_1, \ldots, s_n : \bMM_{g,n} \to \bMM_{g,n+1}$
the sections of $\pi$ arising from the $n$-points of $\bMM_{g,n}$.
Then, $\lambda = \det(\pi_* (\omega_{\bMM_{g,n+1}/\bMM_{g,n}}))$ and
$\psi_i = s_i^*(\omega_{\bMM_{g,n+1}/\bMM_{g,n}})$ for $i=1, \ldots, n$.
Here we set
\begin{align*}
[n] & = \{ 1, \ldots, n\}\quad (\text{note that $[0] = \emptyset$}), \\
\Upsilon_{g,n} & = \{ (i, I) \mid \text{$i \in \ZZ$, $0 \leq i \leq g$
and $I \subseteq [n]$} \} \setminus 
\{ (0, \emptyset), (0, \{1\}), \ldots, (0, \{n\})  \}, \\
\overline{\Upsilon}_{g,n} & = \{ \{(i,I), (j,J)\} \mid 
(i,I), (j,J) \in \Upsilon_{g,n}, i+j=g, I \cap J = \emptyset, I\cup J = [n] \}.
\end{align*}
The boundary
$\Delta = \bMM_{g,n} \setminus \MM_{g,n}$ has 
the following irreducible decomposition:
\[
 \Delta = \Delta_{irr} \cup 
 \bigcup_{\{ (i,I), (j,J) \} \in \overline{\Upsilon}_{g,n}} 
 \Delta_{\{ (i, I), (j, J) \}}.
\]
A general point of $\Delta_{irr}$ represents an $n$-pointed
irreducible stable curve with one node.
A general point of $\Delta_{\{(i,I), (j,J)\}}$ represents an $n$-pointed
stable curve consisting of an $\vert I \vert$-pointed smooth curve 
$C_1$ of genus $i$ and
a $\vert J \vert$-pointed smooth curve $C_2$ of genus $j$ 
meeting transversally
at one point, where $\vert I \vert$-points on $C_1$ 
(resp. $\vert J \vert$-points on $C_2$)
arise from $\{ s_t \}_{t \in I}$ (resp. $\{ s_l \}_{l \in J}$).
Let $\delta_{irr}$ and $\delta_{\{(i,I),(j,J)\}}$ be the classes of
$\Delta_{irr}$ and $\Delta_{\{(i,I),(j,J)\}}$ in 
$\Pic(\bMM_{g,n})_{\QQ}$ respectively.
For our convenience, we denote $\{ (i, I), (g-i, [n]\setminus I) \}$
by $[i,I]$. Moreover, we set
\[
\overline{\Upsilon}^e_{g,n} = \overline{\Upsilon}_{g,n}
\cup \left\{ [0, \{1\}], \ldots, [0, \{n\}] \right\}
\]
and $\delta_{[0, \{i\}]} = -\psi_i$ for $i=1,\ldots, n$.

\section{Comparison of cohomology groups}

In this section, we would like to show the following theorem,
which is crucial for our note.

\begin{Theorem}
\label{thm:rank:neron:severi:alg:stack}
Let $R$ be a discrete valuation ring with $R \subset \CC$, and
$X$ a proper algebraic stack over $R$
\rom{(}see \rom{\S\S\ref{subsec:alg:stack}} for assumptions
of stacks in this note\rom{)}.
We assume that there are \rom{(i)} a finite group $G$,
\rom{(ii)} a smooth, proper and pure dimensional scheme $Y$ over $R$, and
\rom{(iii)} a surjective morphism $\pi : Y \to X$ over $R$ 
with the following properties:
\begin{enumerate}
\renewcommand{\labelenumi}{(\alph{enumi})}
\item $G$ acts on $Y$ over $X$, i.e.,
$\pi \cdot [g] = \pi$ for all $g \in G$.

\item
$X(\CC) \simeq Y(\CC)/G$ as analytic spaces.
\end{enumerate}
Let $\Spec(k) \to \Spec(R)$ be a geometric point of $\Spec(R)$
\rom{(}i.e., $k$ is an algebraically closed field\rom{)}, and
$X_k$ the coarse moduli space
of $X \times_{\Spec(R)} \Spec(k)$.
If $X_k$ is a normal algebraic scheme over $k$, then
\[
\dim_{\QQ_{\ell}} H^i_{et}(X_k, \QQ_{\ell}) \leq \dim_{\QQ} H^i(X(\CC), \QQ)
\]
for every non-negative integer $i$ and
every prime $\ell$ invertible in $k$.
\end{Theorem}

\Proof
We need three lemmas for the proof of the above theorem.

\begin{Lemma}
\label{lem:G:inv:End}
Let $A$ be a commutative ring with the unity, and 
$G$ a finite group such that the order of $G$
is invertible in $A$.
Then, we have the following:
\begin{enumerate}
\renewcommand{\labelenumi}{(\arabic{enumi})}
\item
Let $V$ be an $A$-module such that $G$ acts on $V$ $A$-linearly.
Let $\rho^G_V : V \to V$ be a map given by
$\rho^G_V(x) = (1/\vert G \vert) \sum_{g \in G} [g](x)$. Then, 
$\rho^G_V \in \End_A(V)$ and $\rho^G_V(V) = V^G$.
Moreover, if $B$ is an $A$-algebra and
$G$ acts on $B$ trivially, then
$\rho^G_{V \otimes_A B} = \rho^G_{V} \otimes_A \operatorname{id}_B$.
In particular, $V^G \otimes_A B = (V \otimes_A B)^G$.

\item
Let $f : V \to W$ be a surjective homomorphism
of $A$-modules.
We assume that $G$ acts on $V$ and $W$ $A$-linearly, and
that $f$ is a $G$-homomorphism.
Then, $f(V^G)=W^G$.
\end{enumerate}
\end{Lemma}

\Proof
(1) is obvious.
Let us consider (2).
First of all, since $f$ is a $G$-homomorphism,
we have $f(V^G) \subseteq W^G$. 
Conversely, let us choose
an arbitrary element $w \in W^G$.
Then, there is $v \in V$ with $f(v) = w$ because $f$ is surjective.
Since $w \in W^G$, for each $g \in G$, there
is $x_g \in \Ker(f)$ with $v - [g](v) = x_g$.

Here we claim that
\[
 [h](x_g) = x_{hg} - x_h\quad
\text{for all $g, h \in G$}.
\]
Acting $h$ to the equation $v - [g](v) = x_g$, we have
$[h](v) - [hg](v) = [h](x_g)$. Moreover, $v - [h](v) = x_h$.
Thus,
\[
 x_{hg} = v - [hg](v) = (v - [h](v)) + ([h](v) - [hg](v)) = x_h +
[h](x_g),
\]
which shows us our claim.

We set $x = (1/\vert G \vert) \sum_{g \in G} x_g$. Then, 
for all $h \in G$,
\[
[h](x) = \frac{1}{\vert G \vert} \sum_{g \in G} [h](x_g) =
\frac{1}{\vert G \vert}\sum_{g \in G}(x_{hg} - x_h) = x - x_h.
\]
Thus, if we set $v' = v - x$, then
\[
[h](v') = [h](v) - [h](x) = (v - x_h) - (x - x_h) = v'
\]
for all $h \in G$. Thus, $v' \in V^G$ and $f(v') = f(v) = w$.
Therefore, we can see that $f(V^G) = W^G$.
\QED

\begin{Lemma}
\label{lem:orbifold:cohom:G:inv}
Let $Y$ be a complex manifold, and $G$ a finite group acting on $Y$
holomorphically. Let $X$ be the quotient analytic space $Y/G$ of $Y$ by
the action of $G$, and $\pi : Y \to X$ the canonical morphism. 
Let $\pi^* : H^i(X, \CC) \to H^i(Y, \CC)$ be the homomorphism
of the singular cohomology groups.
Then, for each $i$, $\pi^*$ is injective and its image is 
the $G$-invariant part $H^i(Y, \CC)^G$
of $H^i(Y, \CC)$.
\end{Lemma}

\Proof
Let $A^i(Y)$ be the space of $C^{\infty}$ $i$-forms on $Y$, and
$A^i(Y)^G$ the $G$-invariant part of $A^i(Y)$. Then, it is well known 
(cf. \cite{Kawa}) that
\[
H^i(Y, \CC) = 
\frac{\Ker(d : A^{i}(Y) \to A^{i+1}(Y))}{\Image(d : A^{i-1}(Y) \to A^i(Y))}
\quad\text{and}\quad
H^i(X, \CC) = 
\frac{\Ker(d : A^{i}(Y)^G \to A^{i+1}(Y)^G)}%
{\Image(d : A^{i-1}(Y)^G \to A^i(Y)^G)}.
\]
Note that
\[
\begin{cases}
\Ker(d : A^{i}(Y)^G \to A^{i+1}(Y)^G) =
\Ker(d : A^{i}(Y) \to A^{i+1}(Y)) \cap A^{i}(Y)^G \\
\Image(d : A^{i-1}(Y)^G \to A^i(Y)^G)=
\Image(d : A^{i-1}(Y) \to A^i(Y)) \cap A^{i}(Y)^G.
\end{cases}
\]
In particular, $\pi^*$ is injective.

Let us consider the natural homomorphism
\[
 \alpha : \Ker(d : A^{i}(Y) \to A^{i+1}(Y))
\to H^i(Y, \CC).
\]
By Lemma~\ref{lem:G:inv:End}.2,
\[
\alpha\left(
\Ker(d : A^{i}(Y) \to A^{i+1}(Y)) \cap A^{i}(Y)^G \right)
= H^i(Y, \CC)^G,
\]
which shows us that
$\pi^*(H^i(X, \CC)) = H^i(Y, \CC)^G$.
\QED

\begin{Lemma}
\label{lem:inj:ecoho:finite}
Let $f : Y \to X$ be a finite surjective morphism 
of normal noetherian schemes.
Then, the natural homomorphism
\[
f^* : H_{et}^i(X, \QQ_{\ell}) \to H_{et}^i(Y, \QQ_{\ell})
\]
is injective for every non-negative integer $i$ and
every prime $\ell$ invertible in $H^0(X, \OO_X)$.
\end{Lemma}

\Proof
Clearly, we may assume that $X$ and $Y$ are connected.
Here we claim the following.

\begin{Claim}
\label{claim:lem:inj:ecoho:finite:1}
For every abelian group $\Lambda$, there is a homomorphism
\[
\rho_{f}(\Lambda) : f_*(\Lambda_Y) \to \Lambda_X
\]
with the following properties:
\begin{enumerate}
\renewcommand{\labelenumi}{(\roman{enumi})}
\item
$\rho_{f}(\Lambda) \cdot f^* = \deg(f) \operatorname{id}$,
where $f^* : \Lambda_X \to f_*(\Lambda_Y)$
is the natural injective homomorphism.

\item
Let $\phi : \Lambda' \to \Lambda$ be
a homomorphism of abelian groups. Then, the following diagram is commutative.
\[
\begin{CD}
f_*(\Lambda'_Y) @>{\rho_f(\Lambda')}>> \Lambda'_X \\
@VVV                       @VVV                  \\
f_*(\Lambda_Y) @>{\rho_f(\Lambda)}>> \Lambda_X
\end{CD}
\]
\end{enumerate}
\end{Claim}

Let $K'$ and $K$ be the function fields of $Y$ and $X$ respectively.
First, we assume that $K'$ is separable over $K$.
Let $K''$ be the Galois closure of $K'$ over $K$, and
$G$ the Galois group of $K''/K$.
Moreover, let $\tilde{Y}$ be the normalization of $Y$ in $K''$,
$g : \tilde{Y} \to Y$ the induced morphism, and
$\tilde{f} : \tilde{Y} \overset{g}{\to} Y \overset{f}{\to} X$
the composition of morphisms $g$ and $f$.
We denote by $(K'/K)(K'')$ 
the set of embeddings of $K'$ into $K''$ over $K$; that is,
\[
(K'/K)(K'') = \{ \sigma : K' \hookrightarrow K'' \mid 
\rest{\sigma}{K} = \operatorname{id} \}.
\]
For each $\sigma \in (K'/K)(K'')$,
there is a morphisms $\tilde{\sigma} : \tilde{Y} \to Y$ over $X$ such that
the induced map of function fields is $\sigma$. 
Here let us consider a homomorphism
$\rho' : f_*(\Lambda_Y) \to \tilde{f}_*(\Lambda_{\tilde{Y}})$
given by $\rho'(x) = \sum_{\sigma \in (K'/K)(K'')} \tilde{\sigma}^*(x)$.
It is easy to see that
\[
\Image(\rho') \subseteq \tilde{f}_*(\Lambda_{\tilde{Y}})^G.
\]
Moreover, since $G$ acts transitively on the fibers of $\tilde{Y} \to X$,
we can see
\[
\tilde{f}_*(\Lambda_{\tilde{Y}})^G \subseteq
\Image(\tilde{f}^* : \Lambda_X \to \tilde{f}_*(\Lambda_{\tilde{Y}})).
\]
Thus $\rho'$ gives rise to
a homomorphism
\[
\rho_f(\Lambda) : f_*(\Lambda_Y) \to \Lambda_X.
\]

Next, let us consider a general case.
Let $K_1$ be the separable closure of $K$ in $K'$, and
$Y_1$ the normalization of $X$ in $K_1$. Then, there
are finite morphisms $g : Y \to Y_1$ and $h : Y_1 \to X$
with $f = h \cdot g$. 
Since $g$ is purely inseparable,
$\Lambda_{Y_1} \overset{\sim}{\to} g_*(\Lambda_Y)$.
Thus, $h_*(\Lambda_{Y_1}) \overset{\sim}{\to} f_*(\Lambda_Y)$.
Let $\rho_{\Lambda}(h) : h_*(\Lambda_{Y_1}) \to \Lambda_X$
be a homomorphism as above. Then,
$\rho_{\Lambda}(f)$ is given by $\deg(g) \rho_{\Lambda}(h)$.

The properties (i) and (ii) are obvious by our construction.

\medskip
Let us go back to the proof of our lemma.
Since $f$ is finite, $H^i(Y_{et}, \ZZ/\ell^m\ZZ) = H^i(X_{et}, f_*(\ZZ/\ell^m\ZZ))$.
Thus, by the above claim,
we have a homomorphism $\varrho_m : H^i(Y_{et}, \ZZ/\ell^m\ZZ) \to 
H^i(X_{et}, \ZZ/\ell^m\ZZ)$ such that
$\varrho_m \cdot f^* = \deg(f) \operatorname{id}$.
Here the following diagram is commutative by the property (ii):
\[
\begin{CD}
H^i(Y_{et}, \ZZ/\ell^{m+1}\ZZ) @>{\varrho_{m+1}}>> H^i(X_{et}, \ZZ/\ell^{m+1}\ZZ) \\
@VVV                                             @VVV \\
H^i(Y_{et}, \ZZ/\ell^{m}\ZZ) @>{\varrho_{m}}>> H^i(X_{et}, \ZZ/\ell^{m}\ZZ)
\end{CD}
\]
Thus, we have $\varrho : H_{et}^i(Y, \QQ_{\ell}) \to
H_{et}^i(X, \QQ_{\ell})$ with $\varrho \cdot f^* = \deg(f) \operatorname{id}$.
Therefore, $f^*$ is injective.
\QED

\bigskip
Let us start the proof of Theorem~\ref{thm:rank:neron:severi:alg:stack}.
Let $\bar{\eta}$ be the geometric generic point of $\Spec(R)$,
and $\bar{t}$ the geometric closed point of $\Spec(R)$.
Here we consider two cases:
\begin{enumerate}
\renewcommand{\labelenumi}{(\roman{enumi})}
\item
The image of $\Spec(k) \to \Spec(R)$ is the generic point.

\item
The image of $\Spec(k) \to \Spec(R)$ is the closed point.
\end{enumerate}

\medskip
In the first case,
by the proper base change theorem (cf. \cite[Chapter~I, Theorem~6.1]{FK}),
\[
H^i_{et}(Y_{k}, \QQ_{\ell}) = H^i_{et}(Y_{\bar{\eta}}, \QQ_{\ell}).
\]
Here, by virtue of Lemma~\ref{lem:inj:ecoho:finite},
the natural homomorphism
$H^i_{et}(X_k, \QQ_{\ell}) \to  H^i_{et}(Y_k, \QQ_{\ell})$
is injective and its image is contained in
$H^i_{et}(Y_k, \QQ_{\ell})^G$ because the action of $G$ is given over
$X_k$. Thus,
\[
 \dim_{\QQ_{\ell}} H^i_{et}(X_k, \QQ_{\ell})
 \leq \dim_{\QQ_{\ell}} H^i_{et}(Y_k, \QQ_{\ell})^G =
\dim_{\QQ_{\ell}} H^i_{et}(Y_{\bar{\eta}}, \QQ_{\ell})^G.
\]
Further, by the proper base change theorem,
the comparison theorem (cf. \cite[Chapter~I, Theorem~11.6]{FK}),
Lemma~\ref{lem:G:inv:End}.1 and Lemma~\ref{lem:orbifold:cohom:G:inv},
\begin{align*}
\dim_{\QQ_{\ell}} H^i_{et}(Y_{\bar{\eta}}, \QQ_{\ell})^G & =
\dim_{\QQ_{\ell}} H^i_{et}(Y_{\CC}, \QQ_{\ell})^G = 
\dim_{\QQ_{\ell}} H^i(Y(\CC), \QQ_{\ell})^G \\
&  = \dim_{\CC} H^i(Y(\CC), \CC)^G 
= \dim_{\CC} H^i(X(\CC), \CC).
\end{align*}
Therefore, we get our assertion.

\medskip
In the second case,
by using the proper base change theorem and
Lemma~\ref{lem:inj:ecoho:finite} as before,
we have
\[
 \dim_{\QQ_{\ell}} H^i_{et}(X_k, \QQ_{\ell})
 \leq \dim_{\QQ_{\ell}} H^i_{et}(Y_k, \QQ_{\ell})^G =
\dim_{\QQ_{\ell}} H^i_{et}(Y_{\bar{t}}, \QQ_{\ell})^G
\]
and
\[
\dim_{\QQ_{\ell}} H^i_{et}(Y_{\bar{\eta}}, \QQ_{\ell})^G 
= \dim_{\CC} H^i(X(\CC), \CC).
\]
Therefore, it is sufficient to show that
\renewcommand{\theequation}{\arabic{section}.\arabic{Theorem}}
\addtocounter{Theorem}{1}
\begin{equation}
\label{thm:rank:neron:severi:alg:stack:eqn:1}
\dim_{\QQ_{\ell}} H^i_{et}(Y_{\bar{t}}, \QQ_{\ell})^G =
\dim_{\QQ_{\ell}} H^i_{et}(Y_{\bar{\eta}}, \QQ_{\ell})^G.
\end{equation}
\renewcommand{\theequation}{\arabic{section}.\arabic{Theorem}.\arabic{Claim}}
Indeed, let $f : Y \to \Spec(R)$ be the canonical morphism, and
we set $F = R^i_{et}f_*(\QQ_{\ell})$.
Then, $G$ acts on the sheaf $F$ of \'etale topology.
Namely, for any \'etale neighborhood $U$ of $\Spec(R)$, $G$ acts on $F(U)$, and
for any \'etale morphism $V \to U$ of
\'etale neighborhoods of $\Spec(R)$, the canonical homomorphism
$F(U) \to F(V)$ is a $G$-homomorphism. Thus, the specialization map
\[
s : H^i_{et}(Y_{\bar{t}}, \QQ_{\ell}) = F_{\bar{t}} \to 
F_{\bar{\eta}} =  H^i_{et}(Y_{\bar{\eta}}, \QQ_{\ell})
\]
is a $G$-homomorphism.
On the other hand,
by virtue of the proper-smooth base change theorem
(cf. \cite[Chapter~I, Lemma~8.13]{FK}),
$s$ is bijective.
Thus, we get \eqref{thm:rank:neron:severi:alg:stack:eqn:1},
which completes the proof of Theorem~\ref{thm:rank:neron:severi:alg:stack}.
\QED

\begin{Corollary}
\label{cor:rank:neron:severi:mg}
Let $g$ and $n$ be non-negative integers with $2g-2+n > 0$, and
$\bcMM_{g,n}$ the algebraic stack classifying $n$-pointed
stable curves of genus $g$. 
Then
\[
\dim_{\QQ_{\ell}} H^i_{et}((\bcMM_{g,n})_k, \QQ_{\ell}) \leq
\dim_{\QQ} H^i(\bcMM_{g,n}(\CC), \QQ)
\]
for every algebraically closed field $k$,
every non-negative integer $i$ and every prime $\ell$ invertible in $k$.
\end{Corollary}

\Proof
By virtue of smoothness of
the moduli of curves with non-abelian level structure
due to Looijenga-Pikaart-de Jong-Boggi
(\cite{Lo}, \cite{PJ}, \cite{BP}),
especially by \cite[Proposition~2.6]{BP}, there are
(1) a positive integer $m$,
(2) a finite group $G$,
(3) a smooth, proper and pure dimensional scheme $Y$ over $\ZZ[1/m]$,
and 
(4) a surjective morphism $\pi : Y \to \bcMM_{g,n} \otimes \ZZ[1/m]$
over $\ZZ[1/m]$ such that
(a) $m$ is invertible in $k$,
(b) $G$ acts on $Y$ over $\bcMM_{g,n} \otimes \ZZ[1/m]$
(i.e. $\pi \cdot [g] = \pi$ for all $g \in G$), and that
(c) $\bcMM_{g,n}(\CC) \simeq Y(\CC)/G$ as analytic spaces.
Here, $(\bcMM_{g,n})_k$ is projective. 
Therefore, Theorem~\ref{thm:rank:neron:severi:alg:stack} implies
our corollary.
\QED

\section{Comparison of the $\QQ$-Picard group with the $\QQ$-N\'eron-Severi group}

In this section, we prove the following theorem.

\begin{Theorem}
\label{thm:biject:Pic:NS}
Let $g$ and $n$ be non-negative integers with $2g-2+n > 0$, and
$\bMM_{g,n}$ the moduli space of $n$-pointed
stable curves of genus $g$ over an algebraically closed field $k$.
Then, the natural homomorphism
$\Pic(\bMM_{g,n})_{\QQ} \to \NS(\bMM_{g,n})_{\QQ}$ is bijective.
\end{Theorem}

\Proof
We need to prepare several lemmas.

\begin{Lemma}
\label{lem:normal:finite:Pic:map}
Let $f : Y \to X$ be a finite and surjective
morphism of normal noetherian schemes.
Then, there is a homomorphism $\Nm_{X/Y} : \Pic(Y) \to \Pic(X)$
such that $\Nm_{X/Y}(f^*(L)) = L^{\otimes \deg(f)}$
for all $L \in \Pic(X)$.
\end{Lemma}

\Proof
Let $K'$ and $K$ be the function fields of $Y$ and $X$ respectively.
Let $\Nm : K' \to K$ be the norm map of $K'$ over $K$.
Here we claim that $\Nm : K' \to K$ gives rise to
$\Nm : f_*(\OO_Y^{\times}) \to \OO_X^{\times}$.
This is a local question, so that we may assume that $Y = \Spec(B)$ and
$X = \Spec(A)$. In this case, our assertion means that $\Nm(x) \in A$ for
all $x \in B$. Since $A$ is normal, to see $\Nm(x) \in A$,
it is sufficient to check that $\Nm(x) \in A_P$ for all
$P \in \Spec(A)$ with $\operatorname{ht}(P) = 1$.
Here $B_P$ is flat over $A_P$.
Thus, $B_P$ is free as $A_P$-module. Hence we can see
$\Nm(B_P) \subseteq A_P$. Therefore, we get our claim.

Let $L \in \Pic(Y)$. Then, by \cite[Lecture~10, Lemma~B]{MumLS},
there is an open covering
$\{ U_{\alpha} \}_{\alpha \in I}$ of $X$ such that $\rest{L}{f^{-1}(U_{\alpha})}$ is
a trivial line bundle; i.e.,
there is $\omega_{\alpha} \in L(f^{-1}(U_{\alpha}))$ with
$\rest{L}{f^{-1}(U_{\alpha})} = \OO_{f^{-1}(U_{\alpha})} \omega_{\alpha}$.
Thus, if we set $g_{\alpha\beta} = \omega_{\beta}/\omega_{\alpha}$ for
$\alpha, \beta \in I$, then $g_{\alpha\beta} \in 
\OO_Y^{\times}(f^{-1}(U_{\alpha} \cap U_{\beta}))$, so that
$\Nm(g_{\alpha\beta}) \in \OO_X^{\times}(U_{\alpha} \cap U_{\beta})$.
Therefore, $\{ \Nm(g_{\alpha\beta}) \}$ gives rise to a line bundle
$M$ on $X$. This is the definition of
$\Nm_{X/Y} : \Pic(Y) \to \Pic(X)$. The remaining assertion is obvious
by our construction.
\QED

\begin{Lemma}
\label{lem:Pic:NS:finite:cov}
Let $f : Y \to X$ be a finite and surjective morphism of normal varieties
over an algebraically closed field $k$.
Then, we have the following.
\begin{enumerate}
\renewcommand{\labelenumi}{(\arabic{enumi})}
\item
$f^* : \Pic(X)_{\QQ} \to \Pic(Y)_{\QQ}$ is injective.

\item
If $\Pic(Y)_{\QQ} \to \NS(Y)_{\QQ}$ is injective, then
so is $\Pic(X)_{\QQ} \to \NS(X)_{\QQ}$.
\end{enumerate}
\end{Lemma}

\Proof
(1)
This is a consequence of Lemma~\ref{lem:normal:finite:Pic:map}.

\medskip
(2) Let us consider the following commutative diagram:
\[
\begin{CD}
\Pic(X)_{\QQ} @>>> \NS(X)_{\QQ} \\
@V{f^*}VV          @VV{f^*}V \\
\Pic(Y)_{\QQ} @>>> \NS(Y)_{\QQ}
\end{CD}
\]
By (1), $f^* : \Pic(X)_{\QQ} \to \Pic(Y)_{\QQ}$ is injective.
Thus, we have (2) using the above diagram.
\QED

\begin{Lemma}
\label{lem:weil:cartier:index}
Let $f : Y \to X$ be a finite and surjective morphism of normal 
noetherian schemes.
If $Y$ is locally factorial \rom{(}i.e., $\OO_{Y,y}$
is UFD for all $y \in Y$\rom{)}, then, for any Weil divisor $D$ on $X$,
$\deg(f) D$ is a Cartier divisor.
\end{Lemma}

\Proof
Clearly, we may assume that $D$ is a prime divisor.
Let $D'$ be a Weil divisor associated with the scheme $f^{-1}(D)$.
Then, $D'$ is a Cartier divisor.
Thus, $\Nm_{Y/X}(\OO_Y(D')) \in \Pic(X)$.
Let $X_0$ be a Zariski open set of $X$ such that
$D$ is a Cartier divisor on $X_0$ and $\codim(X \setminus X_0) \geq 2$.
Then, $\Nm_{Y/X}(\OO_Y(D')) = \OO_X(\deg(f)D)$ on $X_0$.
Thus, $\Nm_{Y/X}(\OO_Y(D')) = \OO_X(\deg(f)D)$ on $X$.
In particular, $\OO_X(\deg(f)D)$ is locally free, which means that
$\deg(f)D$ is a Cartier divisor.
\QED

\begin{Lemma}
\label{lem:Pic:NS:isom:sim:connect}
Let $Y$ be a normal projective variety over an algebraically closed field $k$.
If $Y$ is simply connected and
there is a finite and surjective morphism $f : Z \to Y$
of normal projective varieties such that
$Z$ is smooth over $k$, then
the natural homomorphism
\[
\Pic(Y) \otimes \ZZ[1/\deg(f)] \to 
\NS(Y) \otimes \ZZ[1/\deg(f)]
\]
is bijective. Moreover, if $\deg(f)$ is invertible in $k$,
then $\Pic(Y) \to \NS(Y)$ is bijective.
\end{Lemma}

\Proof
We set $n = \deg(f)$ and $P = (\Pic^0_Y)_{red}$, which is
a subgroup scheme of $\Pic_Y$.
For a positive integer $\ell$, 
let $[\ell] : P \to P$ be a homomorphism given by $[\ell](x) = \ell x$.
First we claim the following.

\begin{Claim}
\label{claim:lem:Pic:NS:isom:sim:connect:1}
$[n](P)$ is proper over $k$.
\end{Claim}

For this purpose, it is sufficient to see that
every closed irreducible curve $C$ in $[n](P)$ is proper over $k$.
For a curve $C$ as above, there are
a proper and smooth curve $T$ over $k$,
a non-empty Zariski open set $T_0$ of $T$, and
a morphism $\phi : T_0 \to P$ such that
the Zariski closure of the image
$T_0 \overset{\phi}{\longrightarrow} P \overset{[n]}{\longrightarrow} P$
is $C$.
The morphism $\phi : T_0 \to P \subseteq \Pic_Y$ gives rise to
a line bundle $L_0$ on $T_0 \times Y$ such that
$\phi(x)$ is the class of $\rest{L_0}{\{ x \} \times Y}$ for
all $x \in T_0$.
Let us take a Cartier divisor $D_0$ such that
$\OO_{T_0 \times Y}(D_0) = L_0$.
Then, there is a Weil divisor $D$ on $T \times Y$ with
$\rest{D}{T_0 \times Y} = D_0$.
Here $\operatorname{id}_T \times f : T \times  Z \to T \times Y$
is a finite and surjective morphism of normal varieties such that
$T \times Z$ is smooth over $k$. Note that 
$\deg(\operatorname{id}_T \times f) = \deg(f) = n$.
Thus, by Lemma~\ref{lem:weil:cartier:index}, $nD$ is a Cartier divisor.
Let $\phi' : T \to P$ be a morphism given by a line bundle
$\OO_{T \times Y}(nD)$. Then,
$\rest{\phi'}{T_0} = [n] \cdot \phi$.
Here $C$ is closed in $P$ because $[n](P)$ is closed in $P$.
Thus, we can see that $\phi'(T) = C$.
Hence $C$ is proper over $k$.

\medskip
Next we claim the following.

\begin{Claim}
\label{claim:lem:Pic:NS:isom:sim:connect:2}
$[\ell] : P(k) \to P(k)$ is injective
for every positive integer $\ell$ invertible in $k$.
\end{Claim}

Since $Y$ is simply connected, by \cite[Proposition~2.11]{FK},
$H^1(Y, (\mu_{\ell})_Y) = 0$.
Thus, the Kummer exact sequence:
\[
1 \to (\mu_{\ell})_Y \to \OO_Y^{\times} \overset{\ell}{\longrightarrow} \OO_Y^{\times} \to 1
\]
yields an injection $[\ell] : \Pic(Y) \to \Pic(Y)$.
Hence $[\ell] : P(k) \to P(k)$
is injective.

\medskip
By Claim~\ref{claim:lem:Pic:NS:isom:sim:connect:2},
$[\ell] : [n](P)(k) \to [n](P)(k)$ is injective for every positive integer $\ell$
invertible in $k$.
Thus, $[n](P) = \{ 0 \}$ because $[n](P)$ is an abelian variety
by Claim~\ref{claim:lem:Pic:NS:isom:sim:connect:1}.
Hence
\[
P(k) \otimes \ZZ[1/n] = \{ 0 \}.
\]
Therefore, we get the first assertion because
$\NS(X) \otimes \ZZ[1/n] = \Pic_Y(k) \otimes \ZZ[1/n]/P(k) \otimes \ZZ[1/n]$.

Moreover, if $n$ is invertible in $k$,
then $[n] : P(k) \to P(k)$ is injective.
Therefore, $P(k) = \{ 0 \}$. Hence, we have the second assertion.
\QED

\begin{Lemma}
\label{lem:num:fiber:rational}
Let $f : Y \to X$ be a morphism of
projective normal varieties over an algebraically closed field $k$.
We assume that \rom{(1)} $X$ and $Y$ are $\QQ$-factorial, \rom{(2)}
$\dim f^{-1}(x) = 1$ for all $x \in X$,
and that \rom{(3)} there is a non-empty open set $X_0$ such that
$f^{-1}(x)$ is a smooth rational curve for every $x \in X_0(k)$.
If $D$ is a $\QQ$-divisor on $Y$ with $D \equiv 0$, then
there is a $\QQ$-divisor $E$ on $X$ such that
$f^*(E) \sim_{\QQ} D$ and $E \equiv 0$.
\end{Lemma}

\Proof
Clearly, we may assume that $D$ is a Cartier divisor.
Then, $f_*(\OO_Y(D))$ is a torsion free sheaf of rank $1$
because $\rest{D}{f^{-1}(x)} = \OO_{f^{-1}(x)}$ for all
$x \in X_0(k)$. Thus, there is a divisor $E$ on $X$ such that
$f_*(\OO_Y(D))^{\vee\vee} = \OO_X(E)$.
Considering the natural homomorphism
$f^*f_*(\OO_Y(D)) \to \OO_Y(D)$. We can find an effective divisor
$T$ on $Y$ such that $D \sim_{\QQ} f^*(E) + T$
and $f(T) \subseteq X \setminus X_0$.
Here $(T \cdot C) = 0$ for all curve $C$ with $\dim f(C) = 0$.
Thus, using Zariski's Lemma (cf. Lemma~\ref{lem:int:on:fibers}), 
there are a $\QQ$-divisor
$S$ on $X$ and a Zariski open set $X_1$ of $X$ such that
$f^*(S) \sim_{\QQ} T$ on $f^{-1}(X_1)$ and
$\codim(X \setminus X_1) \geq 2$.
Therefore, $D \sim_{\QQ} f^*(E + S)$ on $f^{-1}(X_1)$.
Here $\codim(Y \setminus f^{-1}(X_1)) \geq 2$. Thus,
$D \sim_{\QQ} f^*(E + S)$ on $Y$.
It is easy to see that $E + S \equiv 0$ using $D \equiv 0$.
\QED

Let us start the proof of Theorem~\ref{thm:biject:Pic:NS}.
First, let us consider the case $g \geq 2$.
Then, by \cite[Proposition~2.6 and Proposition~3.3]{BP}, there are finite and
surjective morphisms
$Z \to Y$ and $Y \to \bMM_{g,n}$ of normal projective varieties over $k$ such that
$Y$ is simply connected and $Z$ is smooth over $k$.
By Lemma~\ref{lem:Pic:NS:finite:cov}.2 and
Lemma~\ref{lem:Pic:NS:isom:sim:connect},
$\Pic(\bMM_{g,n})_{\QQ} \to \NS(\bMM_{g,n})_{\QQ}$ is bijective.

\medskip
Next let us consider the case $g = 0, 1$.
In order to see our assertion, it is sufficient to show that
if $D \equiv 0$ for a $\QQ$-divisor $D$ on $\bMM_{g,n}$, then
$D \sim_{\QQ} 0$. We prove this by induction on $n$.
First, note that $\dim \bMM_{0,3} = 0$, $\bMM_{0,4} = \PP^1_k$ and
$\bMM_{1,1} = \PP^1_k$.
Let us consider $\pi : \bMM_{g,n} \to \bMM_{g,n-1}$.
If $g=0, 1$, then a general fiber of $\pi$ is
a smooth rational curve.
Moreover, $\bMM_{g,n}$ and $\bMM_{g,n-1}$ are $\QQ$-factorial and
$\dim \pi^{-1}(x) = 1$ for all $x \in \bMM_{g,n-1}$.
Thus, by Lemma~\ref{lem:num:fiber:rational}, there is
a $\QQ$-divisor $E$ on $\bMM_{g,n-1}$ such that
$D \sim_{\QQ} \pi^*(E)$ and $E \equiv 0$.
By the hypothesis of induction, we have $E \sim_{\QQ} 0$.
Thus $D \sim_{\QQ} 0$.
\QED

\begin{Corollary}
\label{cor:cycle:map;inj}
Let $g$, $n$, $\bMM_{g,n}$ and $k$ be the same as in
Theorem~\rom{\ref{thm:biject:Pic:NS}}.
Then, the cycle map
\[
\cl^1 : \Pic(\bMM_{g,n}) \otimes \QQ_{\ell} \to H_{et}^2(\bMM_{g,n}, \QQ_{\ell})
\]
is injective for every prime $\ell$ invertible in $k$.
\end{Corollary}

\Proof
Since $\bMM_{g,n}$ is projective,
the kernel of $\NS(\bMM_{g,n}) \to \NSnum(\bMM_{g,n})$ is finite.
Thus, by Theorem~\ref{thm:biject:Pic:NS},
we have
\[
\Pic(\bMM_{g,n})_\QQ \overset{\sim}{\longrightarrow}
\NS(\bMM_{g,n})_\QQ \overset{\sim}{\longrightarrow}
\NSnum(\bMM_{g,n})_\QQ.
\]
Therefore, it is sufficient to show the following lemma.
\QED

\begin{Lemma}
\label{lem:hom:imply:num}
Let $X$ be a proper algebraic spaces over an algebraically closed field $k$, and
$\ell$ a prime invertible in $k$.
Let $\pi : \Pic(X) \otimes \ZZ_{\ell} \to \NS^{\nu}(X) \otimes \ZZ_{\ell}$
be the natural homomorphism and
$\cl^1 : \Pic(X) \otimes \ZZ_{\ell} \to H_{et}^2(X, \ZZ_{\ell})$
the cycle map.
Then, $\Ker(\cl^1) \subseteq \Ker(\pi)$.
In particular,
\[
\Ker(\Pic(X) \otimes \QQ_{\ell} \to H_{et}^2(X, \QQ_{\ell}))\subseteq 
\Ker(\Pic(X) \otimes \QQ_{\ell} \to \NS^{\nu}(X) \otimes \QQ_{\ell}).
\]
\end{Lemma}

\Proof
Let us consider an exact sequence
\[
\Pic(X) \overset{\ell^m}{\longrightarrow} \Pic(X) \to H^2(X_{et}, \ZZ/\ell^m\ZZ)
\]
arising from the Kummer exact sequence
\[
0 \to \ZZ/\ell^m\ZZ \to \OO_X^{\times} \overset{\ell^m}{\longrightarrow} \OO_X^{\times}
\to 0.
\]
Since $\ZZ_{\ell}$ is flat over $\ZZ$, we have an exact sequence
\[
\Pic(X) \otimes \ZZ_{\ell} \overset{\ell^m}{\longrightarrow} 
\Pic(X) \otimes \ZZ_{\ell} \ \overset{\rho_m}{\longrightarrow} 
H^2(X_{et}, \ZZ/\ell^m\ZZ) \otimes \ZZ_{\ell} = H^2(X_{et}, \ZZ/\ell^m\ZZ)
\]
Note that $\cl^1$ is given by
$\varprojlim \rho_m : 
\Pic(X) \otimes \ZZ_{\ell} \to \varprojlim H^2(X_{et}, \ZZ/\ell^m\ZZ)$. 
Thus, if $x \in \Ker(\cl^1)$, then $\rho_m(x) = 0$ for all $m$.
Therefore, there is $y_m \in \Pic(X) \otimes \ZZ_{\ell}$
with $\ell^m y_m = x$. Let $C$ be an irreducible curve on $X$.
Then, $(x \cdot C) = \ell^m(y_m \cdot C)$. Here $(y_m \cdot C) \in \ZZ_{\ell}$.
Thus,
\[
(x \cdot C) \in \bigcap_{m} \ell^m \ZZ_{\ell} = \{ 0 \}.
\]
Thus, $x \in \Ker(\pi)$.
\QED

\section{Linear independence of the tautological classes}

Let $k$ be an algebraically closed field,
$g$ and $n$ non-negative integers with $2g-2+n > 0$, and
$\bMM_{g,n}$ the moduli space of $n$-pointed
stable curves of genus $g$ over $k$.
Then, we have the following
(see \S\S\ref{subsec:tautological:classes} for the definition of 
$\overline{\Upsilon}_{g,n}$, 
$\overline{\Upsilon}^e_{g,n}$ and
the classes $\delta_{\upsilon}$ ($\upsilon \in
\overline{\Upsilon}^e_{g,n}$)):

\begin{Proposition}
\label{prop:linear:ind:NS:g:geq:1}
\begin{enumerate}
\renewcommand{\labelenumi}{(\arabic{enumi})}
\item
If $g \geq 3$, then $\lambda, \delta_{irr}$ and
$\delta_{\upsilon}$'s \rom{(}$\upsilon \in \overline{\Upsilon}^e_{g,n}$\rom{)}
are linearly independent in
$\NS(\bMM_{g,n})_\QQ$.

\item
If $g=2$, then
$\lambda$ and
$\delta_{\upsilon}$'s \rom{(}$\upsilon \in \overline{\Upsilon}^e_{2,n}$\rom{)}
are linearly independent in
$\NS(\bMM_{2,n})_\QQ$.

\item
If $g=1$, then
$\lambda$ and
$\delta_{\upsilon}$'s \rom{(}$\upsilon \in \overline{\Upsilon}_{1,n}$\rom{)}
are linearly independent in
$\NS(\bMM_{1,n})_\QQ$.
\end{enumerate}
\end{Proposition}

\Proof
First of all, by \cite[Theorem~2.2]{GKM}, we have the following.

\begin{enumerate}
\renewcommand{\labelenumi}{(\alph{enumi})}
\item
If $g \geq 3$, then there is
a morphism $\varphi_{irr}:\PP^1_k \to \bMM_{g,n}$
such that
$\deg(\varphi_{irr}^*(\lambda)) = 0$, 
$\deg(\varphi_{irr}^*(\delta_{irr})) = -1$ and
$\deg(\varphi_{irr}(\delta_{\upsilon})) = 0$ ($\forall\  \upsilon \in
\overline{\Upsilon}^e_{g,n}$).

\item
If $g \geq 2$, then, for every $0 \leq i \leq g-2$ and every 
$I \subset \{1,\ldots,n\}$ with $[i,I] \in \overline{\Upsilon}^e_{g,n}$,
there is a morphism $\varphi_{i,I} : \PP^1_k \to \bMM_{g,n}$
such that
$\deg(\varphi_{i,I}^*(\lambda)) = 0$,
$\deg(\varphi_{i,I}^*(\delta_{irr})) = 0$ and
\[
\deg(\varphi_{i,I}^*(\delta_{ \upsilon})) = \begin{cases}
-1 & \text{if $\upsilon = [i,I]$} \\
0 & \text{if $\upsilon \not= [i,I]$}
\end{cases}
\quad(\forall\  \upsilon \in
\overline{\Upsilon}^e_{g,n}).
\]

\item
If $g \geq 2$, then, for every $1 \leq i \leq g-1$ and every 
$I \subset \{1,\ldots,n\}$ with $[i,I] \in \overline{\Upsilon}^e_{g,n}$,
there is a morphism $\varphi_{i,I} : \PP^1_k \to \bMM_{g,n}$
such that
$\deg(\varphi_{i,I}^*(\lambda)) = 0$,
$\deg(\varphi_{i,I}^*(\delta_{irr})) = -2$ and
\[
\deg(\varphi_{i,I}^*(\delta_{\upsilon})) = \begin{cases}
1 & \text{if $\upsilon = [i,I]$} \\
0 & \text{if $\upsilon \not= [i,I]$}
\end{cases}
\quad(\forall\  \upsilon \in
\overline{\Upsilon}^e_{g,n}).
\]

\item
If $g \geq 1$, then, for every $i,j \geq 0$ and every
$I, J \subseteq \{1,\ldots,n\}$ with
$i + j \leq g-1$, $I \cap J = \emptyset$ and 
$[i,I] \in \overline{\Upsilon}^e_{g,n}$,
there is a morphism $\varphi_{i,j,I,J} : \PP^1_k \to \bMM_{g,n}$
such that
$\deg(\varphi_{i,j,I,J}^*(\lambda)) = 0$,
$\deg(\varphi_{i,j,I,J}^*(\delta_{irr})) = 0$ and
\[
\deg(\varphi_{i,j,I,J}^*(\delta_{\upsilon})) = \begin{cases}
1 & \text{if $\upsilon = [i+j,I\cup J]$} \\
-1 & \text{if $\upsilon = [i,I], [j,J]$} \\
0  & \text{otherwise}
\end{cases}
\]
for all $\upsilon \in \overline{\Upsilon}^e_{g,n}$.
\end{enumerate}

\medskip
(1) First, let us consider the case $g \geq 3$.
We assume that
\[
D = a \lambda + b_{irr} \delta_{irr} + 
\sum_{\upsilon \in \overline{\Upsilon}^e_{g,n}}
b_{\upsilon} \delta_{\upsilon} \equiv 0.
\]
Then, since $\deg(\varphi_{irr}^*(D)) = 0$, we have $b_{irr} = 0$.
Here $g \geq 3$. Thus, for every $\upsilon \in \overline{\Upsilon}^e_{g,n}$,
we can find $i,I$ such that
$0 \leq i \leq g-2$, $I \subseteq \{1, \ldots, n\}$ and
$\upsilon = [i, I]$. Thus, by the above (b), we can see $b_{\upsilon} = 0$
for all $\upsilon \in \overline{\Upsilon}^e_{g,n}$.
Hence $D = a \lambda \equiv 0$. Therefore, $a = 0$.

\medskip
(2) Next, let us consider the case $g=2$.
We assume that
\[
D = a \lambda + 
\sum_{\upsilon \in \overline{\Upsilon}^e_{2,n}}
b_{\upsilon} \delta_{\upsilon} \equiv 0.
\]
If $\upsilon = [0, I]$ for some $I \subseteq \{1,\ldots, n\}$, then,
using (b), we can see $b_{\upsilon} = 0$. Otherwise, we can set
$\upsilon = [1, I']$ for some $I' \subseteq \{1,\ldots, n\}$.
Then, $b_{\upsilon} = 0$ by (c).
Thus $D = a \lambda \equiv 0$. Hence, $a = 0$.

\medskip
(3) Finally, let us consider the case $g=1$.
We assume
\[
D = \lambda + 
\sum_{\upsilon \in \overline{\Upsilon}^e_{1,n}}
b_{\upsilon} \delta_{\upsilon} \equiv 0,
\]
where $b_{[0,\{i\}]} = 0$ for all $i = 1, \ldots, n$.
Then, by virtue of (d),
For every non-empty $I, J \subseteq \{1,\ldots,n\}$ with
$I \cap J = \emptyset$,
$b_{[0,I\cup J]} = b_{[0,I]} + b_{[0,J]}$.
Therefore,
\[
b_{[0,I]} = \sum_{i \in I} b_{[0,\{i\}]} = 0
\]
for all non-empty $I \subseteq \{1, \ldots, n\}$.
Hence $D = a \lambda \equiv 0$. Therefore, $a = 0$.
\QED

\section{Generators of the $\QQ$-Picard group and the cycle map}

Let $g$ and $n$ be non-negative integers with $2g-2+n > 0$, and
$\bcMM_{g,n}$ the algebraic stack classifying $n$-pointed
stable curves of genus $g$. For an algebraically
closed field $k$, let $(\bcMM_{g,n})_k$ be the
coarse moduli scheme of $\bcMM_{g,n} \times_{\Spec(\ZZ)} \Spec(k)$.

\begin{Theorem}
\label{thm:generators:NS:Mgn}
$\Pic((\bcMM_{g,n})_k)_\QQ$ is generated by
\[
\text{$\lambda, \psi_1, \ldots, \psi_n, \delta_{irr}$ and
$\delta_{\upsilon}$'s \rom{(}$\upsilon \in \overline{\Upsilon}_{g,n}$\rom{)}}
\]
for any algebraically closed field $k$.
Moreover, the cycle map
\[
\Pic((\bcMM_{g,n})_k) \otimes \QQ_{\ell} \to H_{et}^2((\bcMM_{g,n})_k, \QQ_{\ell})
\]
is bijective for every prime $\ell$ invertible in $k$.

\end{Theorem}

\Proof
By Corollary~\ref{cor:cycle:map;inj},
the cycle map
\[
\Pic((\bcMM_{g,n})_k) \otimes \QQ_{\ell} \to H_{et}^2((\bcMM_{g,n})_k, \QQ_{\ell})
\]
is injective.
Hence, we get
\[
\dim_{\QQ} \Pic((\bcMM_{g,n})_k)_{\QQ} \leq 
\dim_{\QQ_{\ell}}  H_{et}^2((\bcMM_{g,n})_k, \QQ_{\ell}).
\]
Moreover, by Corollary~\ref{cor:rank:neron:severi:mg},
\[
\dim_{\QQ_{\ell}}  H_{et}^2((\bcMM_{g,n})_k, \QQ_{\ell})
\leq \dim_{\QQ} H^2(\bcMM_{g,n}(\CC), \QQ).
\]
Therefore,
\addtocounter{Claim}{1}
\begin{equation}
\label{eqn:thm:generators:NS:Mgn:1}
\dim_{\QQ} \Pic((\bcMM_{g,n})_k)_{\QQ} \leq \dim_{\QQ} H^2(\bcMM_{g,n}(\CC), \QQ)
\end{equation}
and if the equation holds, then the cycle map is bijective.

In the case $g=0$, the assertions of our theorem are well known
(for example, see \cite{KIntM0n}), so that
it is sufficient to show the following (a)--(c) and
$\dim_{\QQ} \Pic((\bcMM_{g,n})_k)_{\QQ} = \dim_{\QQ} H^2(\bcMM_{g,n}(\CC), \QQ)$
for each case.

\begin{enumerate}
\renewcommand{\labelenumi}{(\alph{enumi})}
\item
If $g \geq 3$, then $\lambda, \psi_1, \ldots, \psi_n, \delta_{irr}$ and
$\delta_{\upsilon}$'s \rom{(}$\upsilon \in \overline{\Upsilon}_{g,n}$\rom{)}
form a basis of $\Pic((\bcMM_{g,n})_k)_\QQ$.

\item
If $g=2$, then
$\lambda, , \psi_1, \ldots, \psi_n$ and
$\delta_{\upsilon}$'s \rom{(}$\upsilon \in \overline{\Upsilon}_{2,n}$\rom{)}
form a basis of
$\Pic((\bcMM_{2,n})_k)_\QQ$.

\item
If $g=1$, then
$\lambda$ and
$\delta_{\upsilon}$'s \rom{(}$\upsilon \in \overline{\Upsilon}_{1,n}$\rom{)}
form a basis of
$\Pic((\bcMM_{1,n})_k)_\QQ$.
\end{enumerate}

First of all, it is well known (cf. \cite{AC}) that
$H^2(\bcMM_{g,n}(\CC), \QQ)$ is generated by
\[
\text{$\lambda, \psi_1, \ldots, \psi_n, \delta_{irr}$ and
$\delta_{\upsilon}$'s ($\upsilon \in \overline{\Upsilon}_{g,n}$)}.
\]

\medskip
(a) By Proposition~\ref{prop:linear:ind:NS:g:geq:1}.1,
$\lambda, \psi_1, \ldots, \psi_n, \delta_{irr}$ and
$\delta_{\upsilon}$'s ($\upsilon \in \overline{\Upsilon}_{g,n}$)
are linearly independent in $H^2(\bcMM_{g,n}(\CC), \QQ)$, and
these are also linearly independent in $\Pic((\bcMM_{g,n})_k)_\QQ$.
Thus, by \eqref{eqn:thm:generators:NS:Mgn:1},
\[
\text{$\lambda, \psi_1, \ldots, \psi_n, \delta_{irr}$ and
$\delta_{\upsilon}$'s ($\upsilon \in \overline{\Upsilon}_{g,n}$)}
\]
gives rise to a basis of $\Pic((\bcMM_{g,n})_k)_\QQ$, and
$\dim_{\QQ} \Pic((\bcMM_{g,n})_k)_{\QQ} = \dim_{\QQ} H^2(\bcMM_{g,n}(\CC), \QQ)$.

\medskip
(b) We know that $10\lambda = \delta_{irr} + 2 \delta_1$ on $\bcMM_{2}(\CC)$.
Thus, 
$H^2(\bcMM_{2,n}(\CC), \QQ)$ is generated by
\[
\text{$\lambda, \psi_1, \ldots, \psi_n$ and
$\delta_{\upsilon}$'s ($\upsilon \in \overline{\Upsilon}_{2,n}$)}.
\]
By Proposition~\ref{prop:linear:ind:NS:g:geq:1}.2,
these are linearly independent in $H^2(\bcMM_{2,n}(\CC), \QQ)$.
Moreover,
these are also linearly independent in $\Pic((\bcMM_{2,n})_k)_\QQ$.
Thus, \eqref{eqn:thm:generators:NS:Mgn:1} shows us that
these form a basis of $\NS((\bcMM_{2,n})_k)_\QQ$ and that
$\dim_{\QQ} \Pic((\bcMM_{2,n})_k)_{\QQ} = \dim_{\QQ} H^2(\bcMM_{2,n}(\CC), \QQ)$.

\medskip
(c) First, we know $\delta_{irr} = 12 \lambda$ on $\bcMM_{1,1}(\CC)$ and
$\psi_i = \lambda + \sum_{i \in I, \vert I \vert \geq 2} \delta_{[0,I]}$
for all $i$ on $\bcMM_{1,n}(\CC)$.
Thus,
$H^2(\bcMM_{1,n}(\CC), \QQ)$ is generated by
$\lambda$ and
$\delta_{\upsilon}$'s ($\upsilon \in \overline{\Upsilon}_{1,n}$).
Therefore, by using Proposition~\ref{prop:linear:ind:NS:g:geq:1}.3
and \eqref{eqn:thm:generators:NS:Mgn:1},
$\lambda$ and
$\delta_{\upsilon}$'s ($\upsilon \in \overline{\Upsilon}_{1,n}$)
form a basis of $\Pic((\bcMM_{1,n})_k)_\QQ$ and
$\dim_{\QQ} \Pic((\bcMM_{1,n})_k)_{\QQ} = \dim_{\QQ} H^2(\bcMM_{1,n}(\CC), \QQ)$.
\QED

\begin{Corollary}
\label{cor:generators:NS:Mgn}
Let $g$ and $n$ be non-negative integers with $2g-2+n > 0$, and
$\MM_{g,n}$ the moduli space of $n$-pointed
smooth curves of genus $g$ over an algebraically closed field $k$.
Then, $\Pic(\MM_{g,n})_\QQ$ is generated by
$\lambda, \psi_1, \ldots, \psi_n$.
\end{Corollary}

\Proof
Let us consider the restriction map
$\Pic(\bMM_{g,n})_{\QQ} \to \Pic(\MM_{g,n})_{\QQ}$.
Since $\bMM_{g,n}$ is $\QQ$-factorial, it is surjective.
Thus, Theorem~\ref{thm:generators:NS:Mgn} implies our
corollary.
\QED

\renewcommand{\thesection}{Appendix \Alph{section}}
\renewcommand{\theTheorem}{\Alph{section}.\arabic{Theorem}}
\renewcommand{\theClaim}{\Alph{section}.\arabic{Theorem}.\arabic{Claim}}
\renewcommand{\theequation}{\Alph{section}.\arabic{Theorem}.\arabic{Claim}}
\setcounter{section}{0}

\section{Zariski's lemma for integral scheme}

Let $R$ be a discrete valuation ring, and
$f : Y \to \Spec(R)$ a flat and projective integral scheme over $R$.
Let $\eta$ be the generic point of $\Spec(R)$ and
$o$ the closed point of $\Spec(R)$.
We assume that the genetic fiber $Y_{\eta}$ of $f$ is geometrically
reduced and irreducible curve. Let $Y_o$ be the special fiber of $f$, i.e.,
$Y_o = f^*(o)$.
Let us consider a paring
\[
\Pic(Y) \otimes \Chow^0(Y_o) \to \Chow^{1}(Y_o)
\]
given by the composition of homomorphisms
\[
\Pic(Y) \otimes \Chow^0(Y_o) \to \Pic(Y_o) \otimes \Chow^0(Y_o) \to
\Chow^{1}(Y_o).
\]
We denote by $x \cdot z$ the image of $x \otimes z$ by the above 
homomorphism.
For a Cartier divisor $D$ on $Y$,
the associated cycle of $D$ is denoted by $[D]$, which is an element
of $Z^1(Y)$.
Let us consider the following subgroup $F_c(Y)$ of $Z^0(Y_o)$:
\[
F_c(Y) = \{ x \in Z^0(Y_o) \mid
\text{$x = [D]$ for some Cartier divisor $D$ on $Y$} \}.
\]
For a Cartier divisor $D$ on $Y$ with $[D] \in F_c(Y)$, and
$y \in F_c(Y)$, $D \cdot y$ depend only on $[D]$.
For, if $D'$ is a Cartier divisor on $Y$ with $[D'] = [D]$, and
$E$ is a Cartier divisor on $Y$ with $y = [E]$, then,
by \cite[Theorem~2.4]{Fu},
\[
D \cdot y = E \cdot [D] = E \cdot [D'] = D' \cdot y.
\]
Thus, we can define a bi-linear map
\[
q : F_c(Y) \times F_c(Y) \to \Chow^1(Y_o)
\]
by $q([D], y) = D \cdot y$.
Moreover, \cite[Theorem~2.4]{Fu} says us that $q$ is symmetric; i.e.,
$q(x, y) = q(y, x)$ for all $x, y \in F_c(Y)$.
Then, we define 
the quadratic form $Q$ on $F_c(Y)$ by
\[
Q(x, y) = \deg ( q(x, y) ).
\]
Then, we have the following Zariski's lemma on integral schemes.

\begin{Lemma}[Zariski's lemma for integral scheme]
\label{lem:int:on:fibers}
\begin{enumerate}
\renewcommand{\labelenumi}{(\arabic{enumi})}
\item
$Q([Y_o], x) = 0$ for all $x \in F_c(Y)_{\QQ}$.

\item
$Q(x, x) \leq 0$ for any $x \in F_c(Y)_{\QQ}$.

\item
$Q(x, x) = 0$ if and only if
$x \in \QQ \cdot [Y_o]$.
\end{enumerate}
\end{Lemma}

\Proof
(1):\quad
This is obvious because $\OO_Y(Y_o) \simeq \OO_Y$.

(2) and (3):\quad
If $x \in \QQ \cdot [Y_o]$, then by (1), $Q(x, x) = 0$.
Thus, it is sufficient to prove that
(a) $Q(x, x) \leq 0$ for any $x \in F_c(Y)_{\QQ}$, and that
(b) if $Q(x, x) = 0$, then $x \in \QQ \cdot [Y_o]$.
Here we need the following sublemma.

\begin{Sublemma}
Let $V$ be a finite dimensional vector space over $\RR$, and
$Q$ a quadratic form on $V$. We assume that there are $e \in V$ and
a basis $\{ e_1, \ldots, e_n \}$ of $V$ with the following properties:
\begin{enumerate}
\renewcommand{\labelenumi}{(\roman{enumi})}
\item
If we set $e = a_1 e_1 + \cdots + a_n e_n$, then
$a_i > 0$ for all $i$.

\item
$Q(x, e) \leq 0$ for all $x \in V$.

\item
$Q(e_i, e_j) \geq 0$ for all $i \not= j$.

\item
If we set $S = \{ (i, j) \mid \text{$i \not= j$ and $Q(e_i, e_j) > 0$} \}$,
then, for any $i \not= j$, there is a sequence $i_1, \ldots, i_l$ such that
$i_1 = i$, $i_l = j$, and $(i_t, i_{t+1}) \in S$ for all $1 \leq t < l$.
\end{enumerate}
Then, $Q(x, x) \leq 0$ for all $x \in V$.
Moreover, if $Q(x, x) = 0$ for some $x \not= 0$, then
$x \in \RR e$ and $Q(y, e) = 0$ for all $y \in V$.
\end{Sublemma}

\Proof
Replacing $e_i$ by $a_i e_i$, we may assume that $a_1 = \cdots = a_n = 1$.
If we set $x = x_1 e_1 + \cdots + x_n e_n$,
then, by an easy calculation, we can show
\[
Q(x, x) = \sum_i x_i^2 Q(e_i, e) - \sum_{i < j} (x_i - x_j)^2 Q(e_i, e_j).
\]
Thus, we can easily see our assertions.
\QED

Let us go back to the proof of Lemma~\ref{lem:int:on:fibers}.
First, we assume that $Y$ is regular.
Let $(Y_o)_{red} = E_1 + \cdots + E_n$ be the irreducible
decomposition of $(Y_o)_{red}$. Since $Y$ is regular,
$E_i$'s are Cartier divisors on $Y$ and
$[E_i] \in F_c(Y)$ for all $i$.
Moreover, we can set $Y_o = a_1 E_1 + \cdots + a_n E_n$
for some positive integers $a_1, \ldots, a_n$.
Thus, if we set $e = [Y_o]$ and $e_i = [E_i]$ for $i=1, \ldots, n$,
then (i), (ii) and (iii) in the above sublemma hold.
Moreover, since $Y_o$ is geometrically connected, (iv) also holds.
Thus, we have our assertion in the case where $Y$ is regular.

\medskip
Next, let us consider a general case.
Clearly we may assume that
$x \in F_c(Y)$; i.e., $x = [D]$ for some Cartier divisor $D$ on $Y$.
By virtue of \cite{Lip},
there is a birational morphism $\mu : Y' \to Y$ of
projective schemes over $R$ such that $Y'$ is regular.
Using the projection formula (cf. \cite[(c) of Proposition~2.4]{Fu}),
\[
\deg (\OO_{Y'}(\mu^*(D)) \cdot [\mu^*(D)] ) =
\deg (\OO_Y(D) \cdot [D] ).
\]
Thus, if $Q([\mu^*(D)], [\mu^*(H)]) \leq 0$,
then $Q([D], [D]) \leq 0$.
Moreover, if there is a rational number $\alpha$ such that
$[\mu^*(D)] = \alpha [Y'_o]$, then
$[\mu^*(D)] = \alpha [\mu^*(Y_o)]$. Thus, taking the push-forward $\mu_*$,
we can see that $[D]= \alpha [Y_o]$ in $\Cycle^1(Y)_{\QQ}$.
Hence, we get our lemma.
\QED

\end{document}